\documentclass[12pt,a4paper]{article}
\usepackage[OT1]{fontenc}
\usepackage[latin1]{inputenc}
\usepackage[english]{babel}
\usepackage{amsfonts}
\usepackage{amsthm}
\newtheorem{teo}{Theorem}[section]
\newtheorem{prop}[teo]{Proposition}
\newtheorem{lem}[teo]{Lemma}
\newtheorem{coro}[teo]{Corollary}
\newtheorem{defi}[teo]{Definition}
\theoremstyle{definition}
\newtheorem{rem}[teo]{Remark}

\def\M{{\cal M}}
\def\um{{\cal U}_{\M}}
\def\bh{{\cal B}({\cal H})}
\def\glp{GL(\M,p)}
\def\sp{{\cal{P}}_p}
\def\p{{\cal{P}}}
\def\blp{{\cal B}(L^p)}
\def\ad{{\rm{ad\, }}}
\def\ulp{U(\M,p)}
\def\lp{L^p(\M,\tau)}

\begin{document}

\title{\vspace*{0cm}Spaces of nonpositive curvature arising from a finite algebra\footnote{2000 MSC. Primary 58B20;  Secondary 47C15, 47L07,  54E50.}}


\date{}
\author{Cristian Conde and Gabriel Larotonda}

\maketitle

\abstract{\footnotesize{\noindent  In this paper we introduce a family of examples that can be regarded as spaces of nonpositive curvature, but with the distinct quality that they are not complete as metric spaces. This amounts to the fact that they are modelled on a finite von Neumann  algebra, and the metrics introduced arise from the trace of the algebra. In spite of the noncompleteness of these manifolds, their geometry can be studied from the view-point of metric geometry, and several techniques derived from the functional analysis are applied to gain insight on their geodesic structure.}\footnote{{\bf Keywords and
phrases:} geodesic convexity, non-commutative measure spaces, positive cone, short geodesic, uniform convexity.}}

\section{Introduction}\label{uno}

The view-point adopted to present the examples in this paper is that of the metric geometry, which is well-suited to deal with these infinite-dimensional phenomena. Since Menger and Wald \cite{menger,wald}, who introduced the notions and methods of curves in metric spaces, geodesic length spaces and comparison triangles, there has been several attempts to exploit the intrinsic insight of metric geometry: the foundations of the theory of metric spaces with upper curvature bounds were laid in the 50's with the work of Alexandrov and Busemann \cite{alexandrov,busemann}, but it was not until recently, with the work of  Ballmann, Gromov and Schroeder \cite{ballmann} among others, with their study of \textit{espaces de longueur} of nonpositive curvature, that this subject has shown its true relevance and connections to many areas of modern mathematics, such as operator theory, hyperbolic groups and topology.

Our examples consist of endowing the cone of positive invertible operators of a finite von Neumann algebra with the rectifiable distances induced by the $p$-norms (which are derived from the finite trace of the algebra). The study of the geometry of cones of positive invertible matrices was initiated by Mostow \cite{mostow} in the 50's; his concern was the Riemannian metric arising from the finite trace of the algebra. Later, Corach, Porta and Recht \cite{cpr1,cpr2} studied the geometry of such cones on ${\cal B}{(\cal H})$, the bounded operators on an infinite dimensional Hilbert space ${\cal H}$. The metric they were interested in, is the Finsler metric induced by the uniform norm of ${\cal B}{(\cal H})$ (the term {\em Finsler metric} is used here in a loose sense, since it is not required that the Hessian of the metric is positive definite: moreover, the uniform norm is not even smooth). So the techniques introduced by them are relevant from the view-point of operator theory, but they certainly do not come from the Riemannian geometry. Recently \cite{conde,conde1,la}, we have studied spaces of perturbations of Schatten operators where again, despite the infinite dimensional setting, the metrics introduced come from the (infinite) trace of ${\cal B}{(\cal H})$. The results obtained in those works, together with a factorization theorem for the group of invertible elements of a $C^*$-algebra obtained by Porta and Recht in \cite{pr}, have been extended to Banach-Lie groups in \cite{cl}. A closer relevant precedent is \cite{al} where the (weak) Riemannian case, corresponding to $p=2$, is studied as a geodesic length space. It is worth mentioning that the inner product $\langle v,w\rangle=\tau(vw^*)$ available in that case enabled the introduction of techniques of Riemannian geometry such as comparison triangles; these tools are not available here but, as we have shown in \cite{cl}, we find in the notion of uniform convexity \cite{clarkson} a valuable substitute.

For Banach-Finsler manifolds, the work of Neeb \cite{neeb} gives a natural setting for the study of Cartan-Hadamard manifolds in infinite dimension. He introduced techniques of dissipative operators to describe the condition of nonpositive curvature in its infinitesimal form (i.e. in terms of the differential of the exponential map). These tools can be adapted to the present context, {\em via} operator algebra techniques.

This paper is organized as follows. In Section \ref{back} we introduce the set of positive invertible elements of a finite von Neumann algebra $\M$, together with some (elementary) considerations regarding its smooth manifold structure as a subset of $\M$, and the relevant metrics induced by the trace that we will consider. In Section \ref{convexos} we describe the geodesic structure of such cones of operators with the given Finsler metric, and we characterize convex submanifolds. Finally, in Section \ref{uniconve}, we study the problem of best approximation from a given point to a convex submanifold, via the notion of uniform convexity, and in the process, we establish some inequalities that are a nonlinear variation on the Clarkson-McCarthy's inequalities \cite{clarkson}.

\section{Background}\label{back}

Let $\M$ be a finite von Neumann algebra, and let $\tau$ be a faithful normal tracial state on $\M$. For $1\le p<\infty$, the $p$-metric is given in $\M$ by the trace, $\|x\|_p=\tau(|x|^p)^{\frac{1}{p}}$. Let $L^p=L^p(\M,\tau)$ stand for the the completion of $\M$ relative to the $p$-metric, with the usual identification $L^{\infty}=\M$. We will use $\|\cdot\|$ to denote the usual (uniform) norm of $\M$. Then $L^1$ can be identified with the pre-dual space $\M_*$ of $\M$, the subspace of linear functionals in $\M^*$ which are ultraweakly continuous. For $1<p<\infty$, we have $(L^p)^*=L^q$, where $p^{-1}+q^{-1}=1$. A good reference of the subject is the book of Takesaki \cite{take}. All these identifications are induced by duality via the trace, namely $(v,w) \mapsto\tau(vw^*)$. The $p$-norms are unitarily invariant in the sense that when $u,v\in \um$ (the unitary group of $\M$) then $\|uxv\|_p=\|x\|_p$ for any $x\in \lp$. We will use $\blp$ to denote the set of bounded linear operators \textit{acting} on the Banach space $\lp$, and $B(\M)=B(L^{\infty})$ to indicate the set of bounded linear operators \textit{acting} on $\M$. The involution $*$ of $\M$ extends to an anti-linear isometry $J:\lp\to \lp$, and we will consider the linear space of self-adjoint elements in $\lp$, that is $
L^p_h=\lp_h=\{x\in \lp:Jx=x\}$. We indicate with ${\M}_h$ (resp. $\M_{ah}$) the self-adjoint operators of $\M$ (resp. skew-adjoint), and clearly ${\M}_h=L^p_h\cap \M$.

\subsection{The algebras $B(\M)$ and $\blp$}

In this section we discuss briefly the relations among the different spectra that arise from the various norms considered in this paper.

\begin{rem}\label{spectro}
Let $L_x,R_x:\M\to \M$ stand for the left and right multiplication by $x\in \M$. Then $L_x,R_x\in B(\M)$. The map $L:\M\to B(\M)$ is a faithful representation of $\M$ into a closed Banach subalgebra $L(\M)$ of $B(\M)$. We use $\sigma_{\cal A}(x)$ to indicate the spectrum of the element $x$ relative to the Banach algebra ${\cal A}$. The same remarks hold for $R_x$ instead of $L_x$. Let $\ad x:\M\to \M$ denote the adjoint representation, $\ad x=L_x-R_x$. That is, $\ad x(y)=[x,y]=xy-yx$. Then $\sigma_{B(\M)}(\ad x)\subset \sigma_{\M}(x)-\sigma_{\M}(x)$ because $L_x$ and $R_x$ commute.
\end{rem}

\begin{rem}\label{herm}
Let $\glp$ be the group of invertible elements of $\blp$. It is a Banach-Lie group since it is open there. Let $\ulp$ stand for the group of linear isometries of $\blp$,
$$
\ulp=\{g\in \glp: \|g\|=\|g^{-1}\|\le 1\}.
$$
Here $\|\cdot\|$ denotes the usual supremum norm of operators on a Banach space. Then $\ulp$ is a real Banach-Lie group (not necessarily in the norm topology of $\blp$), with Banach-Lie algebra
$$
Herm({\M},p)=\{T\in \blp: \|e^{s\,T}\|\le 1\; \mbox{ for any } s\in\mathbb R\},
$$
the space of \textit{Hermitian} elements of $\blp$. The key fact here is that if $T\in Herm({\M},p)$ and $s\in \mathbb R$, then by general considerations of the theory of semi-groups and dissipative operators (for instance, see Lemma 3.1 in \cite{lumer}) the operator $1\pm s T$ is invertible and expansive, namely $\|(1\pm s T)z\|_p\ge \|z\|_p$ for any $z\in \lp$.
\end{rem}

\begin{rem}\label{incluido}
If $x\in \M$, the maps $L_x$ and $R_x$ extend to bounded linear operators of $\blp$ with the same norm (less or equal to $\|x\|$), because $\|yxz\|_p\le \|y\|\|x\|_p\|z\|$ whenever $y,z\in \M$.

We shall denote these maps by $\widetilde{L_x}$ and $\widetilde{R_x}$ respectively. Note that the image of $\widetilde{L},\widetilde{R}$ is not necessarily closed in $\blp$. Apparently $\sigma_{\blp}(\widetilde{L_x})\subset \sigma_{\M}(x)$ and $\sigma_{\blp}(\widetilde{R_x})\subset \sigma_{\M}(x)$. Let $\widetilde{\ad}x\in \blp$ be the extension of $\ad x$. Then
$$
\sigma_{\blp}(\widetilde{\ad}x)\subset \sigma_{\blp}(\widetilde{L_x})-\sigma_{\blp}(\widetilde{R_x})\subset\sigma_{\M}(x)-\sigma_{\M}(x).
$$
In particular $\widetilde{\ad}x$ has real spectrum if $x\in {\M}_h$.

Let $\M^{\times}$ stand for the group of invertible elements of $\M$, and $g\in \M^{\times}$. Consider $Ad_g:\M\to \M$ the adjoint action, $Ad_g x=g x g^{-1}$, that is $Ad_g=L_g R_{g^{-1}}$. The identity $Ad_{e^x}=e^{\ad x}$ holds for any $x\in \M$ because $L_x$ and  $R_x$ commute. The same holds for their extensions to $\blp$: when $x\in {\M}$, then
$$
\|e^{i \ad x}z\|_p=\|Ad_{e^{i x}}z\|_p=\|e^{i x}z e^{-i x}\|_p=\|z\|_p
$$
for any $z\in \M$. Then $\|e^{\; i\;\widetilde{\ad}x}z\|_p=\|z\|_p$ for any $z\in L^ p$ because $\M$ is dense there.
\end{rem}

\begin{lem}\label{expa}
If $x\in {\M}_h$, the operator $1+ i\,\,\widetilde{\ad}x\in \blp$ is  expansive and invertible.
\end{lem}
\begin{proof}
By the Remark \ref{incluido}, $ i\;\widetilde{\ad}x\in Herm({\M},p)$. The assertion now follows from Remark \ref{herm}.
\end{proof}

\subsection{Positive invertible elements}

Let $\Pi=\{z\in \mathbb C: Re(z)>0\}$ be the open right half-plane. Let $\Pi_{\M}=\{x\in \M: \sigma_{\M}(x)\subset \Pi\}\subset \M$, which is an open subset of $\M$ since the spectrum map is lower semi-continuous \cite[Theorem 10.20]{rudin}. Let $\p=\Pi_{\M}\cap {\M}_h$ the set of positive invertible elements in $\M$, i.e. the set of elements $a\in \M$ such that $\sigma_{\M}(a)\subset (0,+\infty)$. Clearly $\p$ is an open subset of ${\M}_h$. Since any positive invertible element $a\in \M$ has its spectrum confined to the real interval $(0,+\infty)$, it admits a unique (real analytic) logarithm in ${\M}_h$, that we shall denote $\ln(a)$.

Let us show that the exponential map is a smooth isomorphism of ${\M}_h$ onto $\p$ (with the subspace topology). We first recall a useful (and well-known) expression for the differential of the usual exponential map. If $f:A\to B$ is a smooth map among smooth manifolds, and $TA,TB$ denote the respective tangent bundles, we indicate with $f_*:TA\to TB$ the differential of $f$, and with $f_{*x}$ its specialization at $x\in A$.

\begin{lem}\label{expo}
Let $x,y\in {\M}_h$. Let $\exp(x)=e^x$ be the usual exponential map of $\M$, and let $F$ be the entire function given by $F(z)=z^{-1}\sinh(z)=\sum_{n\ge 0} \frac{z^{2n}}{(2n+1)!}$. Then
$$
(\exp)_{*x}(y)=\int\limits_0^1 e^{(1-t)x}ye^{tx}\, dt,
$$
and $e^{-x/2} (\exp)_{*x}(y)e^{-x/2}= F(\ad x/2)(y).$
\end{lem}

\begin{lem}\label{comun}
Let $\displaystyle F(z)=z^{-1}\sinh(z)$, and let $w\in {\M}_h$. Then $F(\ad w)\in B(\M)$ and it is invertible there. Moreover $F(\ad w)$ is an isomorphism of ${\M}_h$ onto ${\M}_h$.
\end{lem}
\begin{proof}
Clearly $F(\ad w)$ is a bounded map of $\M$ into $\M$. Since $F$ is an entire function, we have $\sigma_{B(\M)}(F(\ad w))=F(\sigma_{B(\M)}(\ad w))$, and since $w\in {\M}_h$, then $\sigma_{B(\M)}(\ad w)\subset \mathbb R$ by Remark \ref{spectro}. Since $F$ maps the real line onto $[1,+\infty)$, then $F(\ad w)$ is invertible in $B(\M)$. The map $F(\ad w)$ sends ${\M}_h$ into ${\M}_h$ because the power series of $F(z)$ involves only real coefficients and even powers of $z$, and $\ad ^2w=\ad w\circ\ad w$ maps ${\M}_h$ into ${\M}_h$. Since $F(\ad w)$ is invertible, it must be an automorphism of ${\M}_h$.
\end{proof}

The following result is also well-known, we include a proof anyway since it is extremely short with the tools at hand.

\begin{prop}\label{inyectiva}
The exponential map $\exp:{\M}_h\to\p$ given by the usual series is a (real analytic) isomorphism onto $\p$, and it has an inverse $\ln:\p\to {\M}_h$, given locally by the usual power series, which is also a (real analytic) isomorphism.
\end{prop}
\begin{proof}
Certainly $\exp$ maps ${\M}_h$ into $\p$ injectively by the well-known properties of the real functional calculus. On the other hand, if $a\in \p$ then $\sigma_{\M}(a)\subset (0,+\infty)$, and then it has a unique real analytic logarithm $\ln(a)\in {\M}_h$ given, for instance, by the Cauchy functional calculus. Then $\exp$ maps ${\M}_h$ onto $\p$. The map $(\exp)_{*x}$ is an isomorphism of ${\M}_h$ onto ${\M}_h$ by Lemma \ref{comun}, since it clearly maps ${\M}_h$ into ${\M}_h$, and it is the composition of the isomorphism $F(\ad x/2)$ with the isomorphism $z\mapsto e^{x/2}\,z\, e^{x/2}$.
\end{proof}

\subsection{The metric spaces $\sp$, $1\le p\le \infty$}

Now $\p$ is a smooth manifold, isomorphic to ${\M}_h$. Consider $\p$
as a subset of $L^p_h$. We make of $\p$ a weak Banach-Finsler $\sp$ manifold by assigning to each tangent space a metric with the $p$-norm:
\begin{equation}\label{metric}
\|x\|_{a,p}=\|a^{-\frac12}x a^{-\frac12}\|_p,\qquad a\in\p,\; x\in {\M}_h.
\end{equation}
\begin{rem}\label{topo}
This metric is continuous in the uniform topology, but note that the tangent spaces are not complete with it (if $p<\infty$). This metric is natural in the sense that it is invariant for the action of the group of invertible elements (see Corollary \ref{inver}). Note that, since $
\|x\|_p\le \|x\|$, then the uniform closure of sets in $\M$ is always contained in the $L^p$ closure.
\end{rem}

We measure the length of rectifiable curves in $\sp$ in the standard fashion, namely
$$
{\ell}_{\sp}(\gamma)=\int_0^1 \|\dot{\gamma}(t)\|_{\gamma(t),p}dt=\int_0^1 \|\gamma(t)^{-\frac12}\dot{\gamma}(t)\gamma(t)^{-\frac12}\|_p \;dt
$$
for any piecewise smooth curve $\gamma:[0,1]\to\sp$. By smooth we mean $C^1$ and with nonzero derivative. We define the geodesic distance between two points $a,b\in\sp$ as the infimum of the lengths of the curves in $\sp$ joining $a$ to $b$,
$$
d_p(a,b)=\inf\{{\ell}_{\sp}(\gamma):\gamma\subset \sp, \gamma \makebox{ is piecewise smooth }, \gamma(0)=a, \; \gamma(1)=b\}.
$$
This defines a semi-finite distance, and for any $a,b\in \sp$  we will exhibit in Section \ref{convexos} a short geodesic joining $a$ to $b$ in $\sp$, so $(\sp,d_p)$ is indeed a metric space. The topology of $\sp$ induced by this metric does not necessarily match the topology of $\sp$ as a subspace of $\M$ ($p<\infty$).

\section{Convex sets and isometries}\label{convexos}

In this section we explore briefly the structure of the smooth convex sets in $\sp$, and as a byproduct we obtain two distinguished classes of isometries. See the survey \cite{eberlein} by P. Eberlein for some background on the group of isometries of (finite dimensional, Riemannian) manifolds of nonpositive curvature.

\begin{defi}
Let $H$ be a closed real subspace of ${\M}_h$. We say that $H$ is a Lie triple system if $[x,[x,y]] =\ad^2 x(y)\in H$ whenever $x,y\in H$. Let $K=\exp(H)\subset \sp$, then we say that $K$ is an exponential set.
\end{defi}

The name \textit{exponential set} has been borrowed from the fundamental paper \cite{pr} by Porta and Recht. It is not hard to see that the condition above is equivalent to $[x,[y,z]]\in H$ for $x,y,z\in H$. In this case we also say that $H$ is \textit{closed under double bracketing}. The following is a well-known result, the proof of Mostow \cite{mostow} for real matrices adapts verbatim to our situation, therefore it is omitted.
\begin{prop}\label{aba}
Let $H$ be a closed real linear subspace of ${\M}_h$ and $K=\exp(H)$. Then $aba\in K$ for any $a,b\in K$ if and only if $H$ is a Lie triple system.
\end{prop}

\begin{rem}\label{curvas}
Let $a,b\in\sp$, then clearly $ a^{-\frac12} b a^{-\frac12}\in \sp$. Then $\gamma_{a,b}:\mathbb R\to \M$, given by
$$
\gamma_{a,b}(t)= a^{\frac12} ( a^{-\frac12} b a^{-\frac12})^t a^{\frac12},
$$
is in fact a smooth curve in $\p$ joining $a=\gamma(0)$ to $b=\gamma(1)$, since $\p =\exp({\M}_h)$. Note that if $a=1$ then $\gamma(t)=b^t=e^{t\ln(b)}$.

If $a,b\in K=\exp(H)$ with $H$ a Lie triple system, then $\gamma_{a,b}(t)\in K$  \textit{for any value of $t\in \mathbb R$}. Note also that
$$
\|\dot{\gamma}(t)\|_{\gamma(t),p} = \|u \ln( a^{-\frac12} b a^{-\frac12}) u^*\|_p =\| \ln( a^{-\frac12} b a^{-\frac12}) \|_p=\|\dot{\gamma}(0)\|_{\gamma(0),p},
$$
where $u=\gamma^{-\frac12}a^{\frac12} ( a^{-\frac12} b a^{-\frac12})^{\frac{t}{2}}\in \um$. Hence, for any $p\ge 1$, the length of these curves is given by the norm of the speed, which is constant and does not depend on $t$:
$$
{\ell}_{\sp}(\gamma_{a,b})=\|\dot{\gamma}_{a,b}(0)\|_{\gamma_{a,b}(0),p}=\| \ln( a^{-\frac12} b a^{-\frac12} ) \|_p,
$$
and also, since $\gamma_{a,b}(1-t)=\gamma_{b,a}(t)$, then ${\ell}_{\sp}(\gamma_{a,b})={\ell}_{\sp}(\gamma_{b,a})$. In particular, if $a=1$, they have length $\|\ln(b)\|_p$.
\end{rem}

\begin{defi}
Let $K\subset\sp$. We say that $K$ is \textit{convex} if, for any given $a,b\in K$, $\gamma_{ab}(t)\in K$ for any $t\in [0,1]$.
\end{defi}

This definition will hold as the natural one, once we prove (in the next section) that these curves are short geodesics. In what follows in this section we proceed as in the (weak) Riemannian case, see \cite{al}.

\begin{defi}
Let us consider, for $a,b\in\sp$, the geodesic symmetries given by $\sigma_a(b)=ab^{-1}a$ for $a,b\in\p$. Then $\sigma_a$ maps $\p$ into $\p$ by Proposition \ref{aba}. Note that since $\sigma_a^2=id$, then $\sigma_a$ is a bijection of $\p$ onto $\p$. These maps are indeed isometries of $\sp$, since for any piecewise smooth curve $\gamma\subset\sp$, then  $\dot{\gamma^{-1}}=-\gamma^{-1}\dot{\gamma}\gamma^{-1}$, hence
$$
\|\frac{\partial}{\partial t}\sigma_a(\gamma)\|_{\sigma_a(\gamma),p}=\| (a^{-1}\gamma a^{-1})^{\frac12}a\gamma^{-1}\dot{\gamma} \gamma^{-1} a (a^{-1}\gamma a^{-1})^{\frac12}\|_p=\| u \gamma^{-\frac12}\dot{\gamma}\gamma^{-\frac12} u^*\|_p,
$$
with $u=(a^{-1}\gamma a^{-1})^{\frac12}a\gamma^{-\frac12}\in \um$. Then $
\|\frac{\partial}{\partial t}\sigma_a(\gamma)\|_{\sigma_a(\gamma),p}=\|\dot{\gamma}\|_{\gamma,p}$.
\end{defi}

\begin{lem}
If $K=\exp(H)$ is an exponential set in $\sp$, with $H$ a Lie triple system, then the geodesic symmetry $\sigma_a:b\mapsto a b^{-1}a$ maps $K$ into $K$ for any $a\in K$.
\end{lem}
\begin{proof}
It is easy to check that $\sigma_a$ maps the curve $\gamma_{a,b}(t)$ onto $\gamma_{a,b}(-t)$ for any $b\in K$. Then $\sigma_a(b)=\gamma_{a,b}(-1)\in K$ by Remark \ref{curvas}.
\end{proof}

\begin{prop}
Let $K=\exp(H)$ with $H$ a closed real linear subspace of ${\M}_h$. Then $K$ is geodesically convex if and only if $H$ has the Lie triple property.
\end{prop}
\begin{proof}
If $H$ is closed under double bracketing, then $K$ is convex by Proposition \ref{aba} and Remark \ref{curvas} above. Let $a,b\in K$, with $K$ convex. Then
$$
aba= a^{\frac32} (  a^{\frac12}b^{-1}a^{\frac12} )^{-1}a^{\frac32}=\sigma_{a^{\frac32}}\circ \sigma_{a^{\frac12}}(b),
$$
which proves that $aba\in K$, and then by Proposition \ref{aba}, $H$ has the Lie triple property.
\end{proof}

From now on we refer to a Lie triple system $H$ as a \textit{LTS}, and to $K=\exp(H)$ as a \textit{convex exponential set}.

Motivated by the definition of convexity we may consider, for $a\in\sp$, a generalization of the usual exponential map. Let $Exp^a:{\M}_h\to\sp$ be given by
$$
Exp^a(x)=a^{\frac12}\exp( a^{-\frac12}x a^{-\frac12}) a^{\frac12}.
$$
Then $Exp^a$ gives a global chart around $a\in \p$; moreover, if $x\in {\M}_h$, then $Exp^a(tx)=\gamma_{a,b}(t)$ where $b=Exp^a(1)$. In Theorem \ref{cortitas} we will show that these curves are short for the geodesic distance, hence the curve $\gamma(t)=Exp^a(tx)$ is a short geodesics starting at $a$ with initial speed $x$. If $b\in \sp$, then taking $x=a^{\frac12}\ln( a^{-\frac12}b a^{-\frac12}) a^{\frac12}\in \M^p_h$ gives $\gamma(1)=b$. Moreover, $\| (Exp^a)_{*x}(y)\|_{a,p}\ge \|y\|_{Exp^a(x),p}$ for any $x,y\in {\M}_h$, a fact that will be proved in Theorem \ref{aumenta} in the next section.

\begin{coro}
Let $K=\exp(H)$, with $H$ a LTS in ${\M}_h$. Then $K\subset\p$ is a manifold with tangent spaces isomorphic to $H$. The maps $Exp^a:a^{\frac12} H a^{\frac12}\to K$ ($a\in K$) are the charts of $K$, and we can identify $T_aK$ with $a^{\frac12}H a^{\frac12}$. If $H$ splits (i.e. ${\M}_h=H\oplus S$ with $S$ a closed supplement of $H$ in ${\M}_h$) then $K$ is a submanifold of $\p$.
\end{coro}
\begin{proof}
Each $Exp^a$ is  a local isomorphism by Lemma \ref{comun}. It is easy to see that these maps are injective, since $\exp$ is injective by Proposition \ref{inyectiva}. Let us show that $Exp^a$ maps $a^{\frac12}Ha^{\frac12}$ onto $K$. If $b\in K$, let $v=a^{\frac12}\ln(a^{-\frac12}ba^{-\frac12})a^{\frac12}$. Then $v\in a^{\frac12}Ha^{\frac12}$, and  $Exp^a(v)=b$. This proves surjectivity. On the other hand, if $v=a^{\frac12}xa^{\frac12}\in a^{\frac12}Ha^{\frac12}$, then $Exp^a(v)=a^{\frac12}e^xa^{\frac12}\in K$ by Proposition \ref{aba}, which shows that $Exp^a$ maps $H$ into $K$. The last assertion follows from the inverse function theorem for Banach spaces, applied to the map $E: H\times S\to \p$ given by $(x,y)\mapsto e^xe^y$.
\end{proof}

\begin{rem}
Let $K=\exp(H)$, with $H$ a LTS. Let $G_K$ be the group generated by the elements in $K$ (namely $g\in G_K$ if $g=a_1\cdots a_n$, with $a_i\in K$). Let $I_g(a)=g^*ag$, for $a\in \p$. Then if $g\in G_K$, $I_g(a)=a_n\cdots a_1 a \, a_1\cdots a_n$ hence $I_g(a)\in K$ by Proposition \ref{aba}. Conversely, every element $a\in K$ can be written as $I_g(1)$, where $g=a^{\frac12}$. Let $U_K=\{g\in G_K:I_g(1)=1\}$ be the isotropy group of $1\in K$ for this action $I$. Then $U_K$ is a subgroup of $\um$ and $K\simeq G_K/U_K$ is an homogeneous space. Moreover, the group $G_K$ acts transitively on $K$, since if $a,b\in K$, then $g=a^{-\frac12}b^{\frac12}\in G_K$ and $I_g(a)=b$.
\end{rem}

\begin{rem}
Let $K=\exp(H)$ with $H$ a LTS. Let $[H,H]\subset \M_{ah}$ denote the closure of the set of finite real linear commutators of elements in $H$. By the Jacobi identity, the real linear space $\mathfrak g=H\oplus [H,H]$ is a real Banach-Lie subalgebra of $\M$. Now $\mathfrak g$ is integrable \cite[Corollary V.2.21]{neeb2} and if $G_H$ denotes the group generated by $\exp(\mathfrak g)$, then $G_H$ is a connected real Banach-Lie group (with a topology that is possibly finer than the norm topology on $\M$), with real Banach-Lie algebra $\mathfrak g$. Note that $G_K\subset G_H$. We claim that $G_H$ also acts on $K$: let us show that $I_g(e^{x_0})=g^*e^{x_0} g$ is an element of $K$ for any $g\in G_H$ and any $x_0\in H$. By the implicit function theorem, any element in $g\in G_H$ can be written as a finite product
$$
g=(e^{x_1}e^{y_1})^{\alpha_1}\cdots (e^{x_n}e^{y_n})^{\alpha_n},
$$
where $x_i\in H$, $y_i\in [H,H]$ and $\alpha_i=\pm 1$. Then direct inspection of the expression  shows that $g^* e^{x_0} g\in K$, because we have either products of the form $e^{x_i}e^{x_j}e^{x_i}$ (which belong to $K$ by Proposition \ref {aba}), or either products of the form $e^{-y_i}e^{x_j}e^{y_i}=\exp(e^{ad(-y_i)}(x_j))$ (which belong to $K$ because $\mathfrak g$ is a Banach-Lie algebra). Let $U_H$ be the unitary part of $G_H$. Then we have the manifold isomorphism $G_H/U_H\simeq K$ given by the action $g\mapsto I_g(1)=g^*g$, that makes $K$ an homogeneous manifold.

These are results on polar decomposition of Banach-Lie groups that are proved in a broader context by K.-H. Neeb \cite{neeb} and these remarks can be obtained by specialization from there.
\end{rem}

\begin{coro}\label{inver}
Let $K=\exp(H)$ be a convex exponential set. Let $g\in G_H$, with $G_H\subset \M^{\times}$ the connected Banach-Lie group with Banach-Lie algebra $\mathfrak g=H\oplus [H,H]$. Let $I_g(a)=g^*ag$. Then each $I_g$ is a bijection of $K$ and an isometry relative to the Finsler metric: the maps $I_g$ act isometrically and transitively on $K$, hence $K$ is an homogeneous Finsler manifold.
\end{coro}
\begin{proof}
Clearly each $I_g$ is a bijection of $K$ by the remarks above. We have $(I_g)_{*a}(x)=g^*xg$ for any $a\in K,x\in H$. Now consider $u=(g^*ag)^{-\frac12} g^* a^{\frac12}$. Then $u\in \um$, hence $\|u y u^*\|_p=\|y \|_p$. Put $y=a^{-\frac12}xa^{-\frac12}$, then $
\|g^*xg\|_{g^*ag,p}=\|x\|_{a,p}$, which proves that the action is isometric. The maps $I_g$ act transitively on $K$, since taking $g=a^{-1/2}b^{1/2}\in G_H$ maps $a$ to $b$.
\end{proof}

\subsection{The exponential metric increasing property}

In this section we prove results related to the existence of short curves for the geodesic distance. Let $\ell_p$ denote the usual $p$-length of curves in the linear space ${\M}_h$.

The considerations of this section are an extension of the results in \cite{neeb}, where the author considers Banach-Finsler manifolds with spray. What is remarkable is that those considerations still hold in this setting, disregarding the fact that the topology of $\sp$ with the $p$-norms does not match the topology of $\p$ as a Banach manifold (if $p<\infty$). The following lemma is Proposition 3.15 in \cite{neeb}, adapted to our situation.

\begin{lem}\label{expande}
Let $F(z)=z^{-1}\sinh(z)$, $w\in {\M}_h$. Then $F(\ad w)$  admits a bounded extension to $\blp$, given by the analytic functional calculus of $\widetilde{\ad}w $, which is invertible and expansive, $\| F(\widetilde{\ad}w)(z)\|_p\ge \|z\|_p$ for any $z\in \lp$.
\end{lem}
\begin{proof}
Let us write $F(z)$ in its Weierstrass expansion. Since the zero set of $F$ is $\{z_k=k\pi i\}$, then
$F(z)=\prod_{n\ge 1}( 1\pm i\frac{z}{n\pi})$, where the product converges uniformly on compact sets of $\mathbb C$ to $F$. Let $T_n= 1\pm i\frac{\widetilde{\ad} w}{n\pi}$, then $T_n\in \blp$ and it is expansive by Lemma \ref{expa}. Hence
$$
F(\widetilde{\ad}w)=\lim\limits_n\prod\limits_{k=1}^n T_k\in \blp
$$
and it is expansive there since each term is expansive.
\end{proof}

\begin{rem}
The inequality of the previous lemma is equivalent to the so-called {\em exponential metric increasing property} \cite{bhatia2}, which states that
$$
\|\int_0^1 a^{1-t}b a^t\, dt\|_p\ge \|a^{\frac12} b a^{\frac12}\|_p.
$$
Indeed, put $b=a^{-\frac12} y a^{-\frac12}$ in the above equation, put $a=e^{x}$, and use the identities of Lemma \ref{expo}.
\end{rem}

Let $\gamma$ be a piecewise smooth curve  $\gamma\subset\p$. Then $\gamma=e^{\Gamma}$ for a uniquely determined piecewise smooth curve $\Gamma=\ln(\gamma)$ such that $\Gamma\subset {\M}_h$. By Lemma \ref{expo},
$$
\dot{\gamma}=(\exp)_{*\Gamma}(\dot{\Gamma})=\int_0^1 e^{(1-t)\Gamma}\dot{\Gamma}e^{t\Gamma}\,dt.
$$

\begin{teo}\label{aumenta}
Let $\gamma=e^{\Gamma}\subset \sp$ be a piecewise smooth curve. Then $\ell_p(\Gamma)\le {\ell}_{\sp}(\gamma)$.
\end{teo}
\begin{proof}
Let us compute the speed of $\gamma$ using Lemma \ref{expo}:
$$
\|\dot{\gamma}\|_{\gamma,p} =\|\gamma^{-1/2}\dot{\gamma}\gamma^{-1/2}\|_p = \|e^{-\frac{\Gamma}{2}} (\exp)_{*\Gamma}(\dot{\Gamma}) e^{-\frac{\Gamma}{2}} \|_p =\| F(\ad \,\Gamma/2)(\dot{\Gamma})\|_p.
$$
On the other hand, by Lemma \ref{expande}, $\| F(\ad \,\Gamma/2)(\dot{\Gamma})\|_p\ge \|\dot{\Gamma}\|_p$.
\end{proof}

\begin{coro}\label{EMI}
Let $v,w\in {\M}_h$. Then $d_p(e^w,e^v)\ge \|w-v\|_p$.
\end{coro}
\begin{proof}
Let $\gamma$ be any piecewise smooth curve joining $e^v$ to $e^w$, put $\gamma=e^{\Gamma}$ as before. Then
$$
\|w-v\|_p=\|\Gamma(1)-\Gamma(0)\|_p=\|\int_0^1\dot{\Gamma}\,dt\|_p
\le \int_0^1 \|\dot{\Gamma}\|_p\, dt = \ell_p(\Gamma)
\le {\ell}_{\sp}(\gamma).
$$
Hence the infimum of the length of these curves must be greater or equal than $\|w-v\|_p$.
\end{proof}

\begin{teo}\label{cortitas}
Let $a,b\in\sp$. Let $\gamma_{a,b}(t)= a^{\frac12} ( a^{-\frac12} b a^{-\frac12})^t a^{\frac12}$. Then $\gamma_{a,b}$ is shorter than any other piecewise smooth curve joining $a$ to $b$ in $\sp$, and
$$
d_p(e^v,e^w)=\|\ln(e^{v/2}e^{-w}e^{v/2})\|_p.
$$
\end{teo}
\begin{proof}
By Corollary \ref{inver}, it suffices to prove the result for $a=1$, and $b=e^w$. Let $\gamma(t)=e^{tw}$. Let $\beta$ be any other curve joining $1$ to $b$. By Corollary \ref{EMI},
$$
\ell_{\sp}(\gamma)=\|w\|_p=\|w-0\|_p\le d_p(e^w,1)\le \ell_{\sp}(\beta).\qedhere
$$
\end{proof}

\begin{rem}
For $1<p<\infty$ the strict convexity properties of $L^p$ imply that straight segments are the unique short smooth curves joining two vectors in $L^p$, and this in turn implies that, for $1<p<\infty$, the geodesics of Theorem \ref{cortitas} are the unique short curves joining $a,b\in \sp$, when the length is measured with the tangent $p$-norms. Indeed, if $\gamma$ is a short smooth curve joining $1$ to $a=e^v$ in $\sp$, then $\gamma=e^{\Gamma}$ for some smooth curve $\Gamma\in {\M}_h$, and since $\ell_p(\Gamma)\le\ell_{\sp}(\gamma)=\|v\|_p$, then $\Gamma(t)=tv$ hence $\gamma=e^{tv}=\gamma_{1,a}(t)$. By the invariance of the metric, the claim follows for $a,b\in \sp$.
\end{rem}

\begin{coro}
Let $K\subset \sp$ be convex. Then for any  $a,b\in K$, the curve $\gamma_{a,b}$ is a smooth short path joining $a$ to $b$ in $K$. It is unique if $1<p<\infty$.
\end{coro}

\begin{rem}
The minimal curves in $\p_p$ have the following property: if we apply the complex interpolation method introduced by Calderón in  \cite{calderon}, to the space $L^p(\M,\tau)$ with the Finsler norms $\|\;\|_{a, p}$ and $\|\;\|_{b, p}$ with $a,b \in \p$, it can be proved that the interpolated curve matches the minimal curve $\gamma_{a,b}$ in $\p$. This is an extension of the result obtained in \cite{interpole} for the Finsler metric induced by the uniform norm in the cone of positive invertible operators of a $C^*$-algebra, but with a different approach. Following the notation used in \cite{bl}, observe that for all $a,b \in \p$ and $1\leq p < \infty$, the Banach spaces $L^p_a=(L^p,\|\;\|_{a,p})$ and  $L^p_b=(L^p,\|\;\|_{b,p})$ are compatible, due to the isomorphism $L^p_a \simeq (L^p,\|\;\|_p)$. Then the following theorem can be proved as in \cite[Theorem 3.1]{condeinterpola}, with some minor adaptations to the proof that therefore, is omitted here:
\begin{teo}
Let $p\ge 1$, $a,b \in \p$ and $t\in (0,1)$. Then $
(L^p_a,L^p_b)_{[t]}=L^p_{\gamma_{a,b}(t)}$.
\end{teo}
\end{rem}

\subsection{Completion of $\p$}

There are three natural metrics to consider in the manifold $\p$. One is the linear metric induced by the $p$-norms, when one regards $\p$ as a linear subspace of $L^p_{h}$, that is
\begin{equation}\label{inducida}
d^l_p(e^v,e^w)=\|e^v-e^w\|_p,
\end{equation}
for $v,w\in \M_h$. The second one is the Finsler metric induced by the $p$-length functional on rectifiable arcs, that is
$$
d_p(e^v,e^w)=\|\ln(e^{v/2}e^{-w}e^{v/2})\|_p.
$$
The third one is the metric induced by the isomorphism of $\p$ with its tangent space $\M_h$, that is
\begin{equation}\label{usual}
d^t_p(e^v,e^w)=\|v-w\|_p.
\end{equation}
By the exponential metric increasing property, one can compare $d_p\ge d^t_p$.

Assume that $v,w\in \M_{h}$, and consider $\alpha(t)=\ln(e^{tv/2}e^{-tw}e^{tv/2})$, which is a smooth curve of self-adjoint elements of $\M$. Note that $\alpha(1)=\ln(e^{v/2}e^{-w}e^{v/2})$ and $\alpha(0)=0$. Let $\beta=e^{\alpha}$, then by Theorem \ref{aumenta} followed by Hölder's inequality,
\begin{eqnarray}
d_p(e^v,e^w)&= & \|\alpha(1)-\alpha(0)\|_p\le \int_0^1 \|\dot{\alpha}(t)\|_p dt \le \int_0^1 \|\beta^{-1/2}\dot{\beta}\beta^{-1/2}\|_p dt \nonumber\\
& \le &  \int_0^1 \|\beta^{-1}\dot{\beta}\|_p dt \le \int_0^1 \|\beta^{-1}\| \| \dot{\beta}\|_p dt\nonumber
\end{eqnarray}
We have used the elementary inequality $\|xy\|_p\le \|yx\|_p$ if $(xy)^*=xy$.
A straightforward computation shows that
$$
\dot{\beta}(t)=\frac12 e^{tv/2}(v-w)e^{-tw}e^{tv/2}+\frac12 e^{tv/2}e^{-tw}(v-w)e^{tv/2}
$$
Thus $\|\dot{\beta}(t)\|_p\le e^{t(\|v\|+\|w\|)} \|v-w\|_p$, since $
\|e^{tw}\|\le e^{t\|w\|}$. Likewise,
$$
\|\beta^{-1}\|=\|e^{-tv/2}e^{tw}e^{-tv/2}\|\le e^{t(\|v\|+\|w\|)}.
$$
It follows that
$$
d_p(e^v,e^w)\le \int_0^1 e^{2t(\|v\|+\|w\|)}dt\; \|v-w\|_p\le \frac{e^{2(\|v\|+\|w\|)}-1}{2(\|v\|+\|w\|)}\, \|v-w\|_p,
$$
namely
$$
\|v-w\|_p\le d_p(e^v,e^w)\le K_{\infty}(v,w) \|v-w\|_p
$$
where $K_{\infty}$ is a constant depending solely on the uniform norm of $v,w$, such that $K_{\infty}(v,w)\to 1$ as $v,w\to 0$ in $\M$.

Now we compare $\|e^v-e^w\|$ to $\|v-w\|_p$. Let $v,w\in \M_h$, with $\|v\|,\|w\|\le C$. Then
\begin{eqnarray}
\|e^v-e^w\|_p &= &\|\sum\limits_{n\ge 1} \frac{v^n}{n!}-\frac{w^n}{n!}\|_p=\|\sum\limits_{n\ge 1}\frac{1}{n!} \sum\limits_{j=0}^{n-1} v^{n-1-j}(v-w)w^j\|_p\nonumber\\
&\le & \sum\limits_{n\ge 1}\frac{1}{n!} \sum\limits_{j=0}^{n-1} \|v\|^{n-1-j}\|v-w\|_p\|w\|^j\le \sum\limits_{n\ge 1}\frac{1}{(n-1)!}C^{n-1} \|v-w\|_p\nonumber\\
&=&e^C \|v-w\|_p.\nonumber
\end{eqnarray}
Likewise, if $m=\max\{\|e^v\|,\|e^w\|\}$ and $\delta=\max\{\|\frac{e^v}{m}-1\|,\|\frac{e^w}{m}-1\|\}<1$, then the expansion $\ln(x)=-\sum\limits_{n\ge 1} \frac1n (x-1)^n$ for any $x$ such that $|x-1|<1$ gives $\|v-w\|_p\le \frac{1}{m(1-\delta)}\|e^v-e^w\|_p$. Thus
$$
d^t_p\le C\, d^l_p\le C'\,  d^t_p\le C'\, d_p\le C''\, d^t_p.
$$
for uniformly bounded subsets of $\p$, where the three metrics are equivalent. Then such subsets of $\p$ are complete with the distance $d_p$: they are complete with the linear $p$-metric $d^l_p(v,w)=\|v-w\|_p$ in $\M$, since the linear $p$-metric induces the strong operator topology on uniformly bounded subsets of $\M$.


It is not hard to see that the completion of $\p$ with the distance (\ref{inducida}) gives the positive (non-necesarily invertible, possibly unbounded) operators of $L^p(\M,\tau)$.

\begin{rem}
What is not clear, and we would like to know, is the structure of the completion of $\p_p$ relative to its rectifiable metric. In the finite dimensional setting, it is well-known that $\p_p$ is a complete metric space with it. Certainly, it is not a complete space: for consider $x\in L^p$ such that $Jx=x$ and $x$ is an unbounded operator affiliated with $\M$. For $n\in\mathbb N$, let $p_n$ be the spectral projection of $x$ obtained from the finite interval $[-n,n]\subset\mathbb R$, and let $x_n=xp_n=p_nx$. Then it is easy to check that $x_n\in {\cal M}$, $x_n^*=x_n$, $x_n x_m=x_m x_n$ for any $n,m \in\mathbb N$ and moreover $\|x_n-x\|_p\to 0$. Thus $d_p(e^{x_n},e^{x_m})=\|x_n-x_m\|_p<\epsilon$ if $n,m>n_0$, but $e^{x_n}$ cannot converge to a point in $\p$.
\end{rem}

\subsection{Cördes inequality and convexity of the geodesic distance}

The Cördes inequality for bounded operators on $\bh$  states that $
\|e^{tx}e^{ty}\|\le \|e^xe^y\|^t$ for any self-adjoint $x,y$ and  $t\in [0,1]$. It is equivalent to the inequality
$$
\|\ln(e^{tx}e^{-2ty}e^{tx})\|\le t \|\ln(e^{x}e^{-2y}e^{x})\|,
$$
which has a geometric interpretation \cite{cpr}: it establishes the fact that the geodesic distance in the space of positive invertible operators is a convex function. In our context it can be related to a well-known inequality due to Araki, Lieb and Thirring \cite{kosa}.

\begin{lem}\label{unmedio}
Let $x,y\in {\M}_h$. Then for any $t\in [0,1]$,
$$
\|\ln(e^{\frac{tx}{2}}e^{-ty} e^{\frac{tx}{2}})\|_p \le t\; \|\ln(e^{\frac{x}{2}}e^{-y}e^{\frac{x}{2}})\|_p.
$$
\end{lem}
\begin{proof}
See for instance \cite{inelar} for a detailed proof.
\end{proof}

The following fact was proved in \cite{lim} by Lawson and Lim for Banach-Finsler manifolds, it still holds in our weak setting due to the previous lemma.

\begin{teo}
Let $\gamma_{a,b}$ and $\gamma_{a,c}$ be two short curves as in Remark \ref{curvas}, starting both at $a\in\sp$. Let $f(t)=d_p(\gamma_{a,b}(t),\gamma_{a,c}(t))$ be the distance function among the two geodesics. Then  $f:[0,1]\to \mathbb R_{\ge 0}$ is continuous and convex.
\end{teo}
\begin{proof}
By the invariance of the metric, it will suffice to prove the theorem assuming $a=1$, where
$f(t)=\|\ln(c^{\frac{t}{2}}b^{-t}c^{\frac{t}{2}})\|_p$. It is continuous, since $t\mapsto h(t)=\ln(c^{\frac{t}{2}}b^{-t}c^{\frac{t}{2}})$ is continuous as a map from $[0,1]$ to $\M$ with the uniform topology, and
$$
|f(s)-f(t)|\le\|h(s)-h(t)\|_p\le \|h(s)-h(t)\|.
$$
The convexity of $f$ is equivalent to the inequality
$$
f(t)=\|\ln(c^{\frac{t}{2}}b^{-t}c^{\frac{t}{2}})\|_p\le t \|\ln(c^{\frac12}b^{-1}c^{\frac12})\|_p=tf(1)
$$
for any $t\in (0,1)$, which is exactly the claim of Lemma \ref{unmedio}.
\end{proof}

\begin{coro} Let $\gamma_{a,b}$ and $\gamma_{c,d}$ be two short curves in $\sp$ as in Remark \ref{curvas}. Let $g:[0,1]\to\mathbb R_{\ge 0}$ be the distance among the two geodesics. Then $g$ is continuous and convex.
\end{coro}
\begin{proof}
The map $g$ is continuous by the same argument we used in the proof of the previous theorem. Now consider the geodesic rectangle with vertices $a,b,c,d$, let $\gamma_{c,b}$ be the short curve joining $c$ to $b$ in $\sp$, and consider the triangle with sides $c,b,d$ and the geodesic triangle with sides $b,a,c$. Note that $\gamma_{c,b}(t)=\gamma_{b,c}(1-t)$ and the same holds for $\gamma_{a,b}$. Then, by the triangle inequality
$$
g(t)=d_p(\gamma_{a,b}(t),\gamma_{c,d}(t))\le d_p(\gamma_{a,b}(t),\gamma_{c,b}(t))+d_p(\gamma_{c,b}(t),\gamma_{c,d}(t)).
$$
By the previous theorem $d_p(\gamma_{c,b}(t),\gamma_{c,d}(t))\le t\, d_p(b,d)$, and also
$$d_p(\gamma_{b,c}(1-t),\gamma_{b,a}(1-t))\le (1-t)\, d_p(a,b).
$$
Adding these two inequalities yields the convexity of $g$.
\end{proof}

\begin{coro}\label{puntocurva}
Let $a\in\sp$, and $\gamma_{b,c}$ a short curve as in Remark \ref{curvas}. Then the distance map $f:[0,1]\to \mathbb R_{\ge 0}$ from the point $a$ to the curve $\gamma_{b,c}$ is a continuous and convex function.
\end{coro}

\section{Uniform convexity and minimizers}\label{uniconve}

The notion of uniform convexity for Banach spaces was introduced in \cite{clarkson}, where Clarkson showed that the classical measure spaces $L^p(\Omega,\mu)$ ($1<p<\infty)$ are uniformly convex. This notion can be translated to inner metric spaces. The notion of midpoint map plays a fundamental role. Let $(X,d)$ be a metric space. A \textit{midpoint map} $m:X\times X \to X$ is an assignment satisfying
$$
d(m(x,y),x)=\frac 12 d(x,y)=d(m(x,y),y)\qquad \forall x,y \in X.
$$

\begin{defi}Let $(X,d)$ be a metric space with a midpoint map $m$. Then $X$ is uniformly ball convex if for all $0<\epsilon\le 2$ there exists  $\delta_d(\epsilon)>0$ such that for all $x,y,z\in X$ satisfying $d(x,y)>\epsilon\; max\{d(x,z),d(y,z)\}$, it holds
$$
d(m(x,y),z)\leq(1-\delta_d(\epsilon)) {\rm{max}}\{d(x,z),d(y,z)\}.
$$
The function $\delta_d$ is called the modulus of convexity of the space.
\end{defi}

\subsection{Uniform convexity of  $\sp$, with $1<p<\infty$.}

Clarkson-McCarthy's inequalities do hold in non commutative $L^p$ spaces, as shown by Kosaki, see \cite[Propositions 5.2 and 5.3]{kosaki}:
\begin{prop}\label{kosa}
Let $a,b\in L^p(M,\tau)$, $1< p\le 2$, and $1/p+1/q=1$. Then
$$
(\|a+b\|_p^q+\|a-b\|_p^q)^{\frac1q}\le 2^{\frac1q} (\|a\|_p^p+\|b\|_p^p)^{\frac1p}.
$$
If $2\le p<\infty$ then
$$
(\|a+b\|_p^p+\|a-b\|_p^p)^{\frac1p}\le 2^{\frac1q} (\|a\|_p^p+\|b\|_p^p)^{\frac1p}.
$$
\end{prop}

\begin{lem}\label{pmenor}
For $x,y \in L^p(\M,\tau)$, $1<p\le 2$, and $1/p+1/q=1$, we have
$$
\left\|x\right\|_p^q+\left\|y\right\|_p^q\leq \frac{1}{2}(\left\|x+y\right\|_p^q+\left\|x-y\right\|_p^q).
$$
If $2\le p\le \infty$ then
$$
\left\|x\right\|_p^p+\left\|y\right\|_p^p\leq \frac{1}{2}(\left\|x+y\right\|_p^p+\left\|x-y\right\|_p^p).
$$
\end{lem}
\begin{proof}First we consider $p\le 2$. From the previous proposition, by setting $a=\frac{x+y}{2}$ and $b=\frac{x-y}{2}$ we obtain that
$$
\left\|x\right\|_p^q+\left\|y\right\|_p^q \leq 2(\frac{1}{2^p})^{q/p}(\left\|x+y\right\|_p^p+\left\|x-y\right\|_p^p)^{q/p}=2^{1-q}(\left\|x+y\right\|_p^p+\left\|x-y\right\|_p^p)^{q/p}.
$$
Since $q-1=q/p$, and using the fact that the function $f(t)=t^{q/p}$ is convex on $[0,+\infty)$,
$$
\left\|x\right\|_p^q+\left\|y\right\|_p^q\leq \frac{1}{2}(\left\|x+y\right\|_p^p)^{q/p}+\frac{1}{2}(\left\|x-y\right\|_p^p)^{q/p}=\frac{1}{2}\left\|x+y\right\|_p^q+\frac{1}{2}\left\|x-y\right\|_p^q.
$$

The case $p\ge 2$ is easier, just put $a=\frac{x+y}{2}$, $b=\frac{x-y}{2}$ in Proposition \ref{kosa} to obtain it.
\end{proof}

The following inequalities establish semi-parallelogram laws in $\sp$. We use $\gamma_t$ to indicate the point $\gamma(t) \in\sp$, and $q\ge 1$ denotes the conjugate exponent of $p$ as before, $1/p+1/q=1$. In what follows, $r\ge 2$ indicates the following number: $r=q$ if $p\in (1,2]$ and $r=p$ if $p\in [2,+\infty)$. That is $r=\max\{p,q\}.$

\begin{teo}\label{conveunmedio}
Let $a\in \sp$ and $\gamma:[0,1]\rightarrow \sp$ be a geodesic. Then if $1<p\leq2$,
$$
\frac{1}{2^q}d_p(\gamma_0,\gamma_1)^q \leq \frac{1}{2}(d_p(a,\gamma_0)^q+d_p(a,\gamma_1)^q)-d_p(a,\gamma_{1/2})^q.
$$
If $2\leq p<\infty$ then
$$
\frac{1}{2^p}d_p(\gamma_0,\gamma_1)^p \leq \frac{1}{2}(d_p(a,\gamma_0)^p+d_p(a,\gamma_1)^p)-d_p(a,\gamma_{1/2})^p.
$$
\end{teo}
\begin{proof}Consider $2\leq p<\infty$. By the invariance of the metric, it suffices to consider the case $\gamma_{1/2}=1$. Let $\gamma_0=e^y$, $\gamma_1=e^{-y}$, $a=e^x$. Then $d_p(\gamma_0,\gamma_1)=\ell_{\sp}(\gamma)=\|2y\|_p=2\|y\|_p$, and $\|x\|_p=d_p(\gamma_{1/2},a)$. By the previous lemma and the exponential metric increasing property of Corollary \ref{EMI},
\begin{eqnarray}
\frac{1}{2^p}d_p(\gamma_0,\gamma_1)^p+ d_p(\gamma_{1/2},e^x)^p=\|y\|_p^p+\|x\|_p^p &\le & \frac12 (\|x +y\|_p^p +\|x-y\|_p^p)\nonumber\\
&\le &\frac12 ( d_p(e^x,e^y)^p + d_p(e^{-y},e^x)^p)\nonumber\\
& = &\frac12 ( d_p(e^x,\gamma_0)^p + d_p(e^x,\gamma_1)^p).\nonumber
\end{eqnarray}
The other inequality has an analogous proof  and it is therefore omitted.
\end{proof}

\begin{coro}\label{bco}
For $1<p<\infty$ the metric space $(\sp,d_p)$ is uniformly ball convex, and an admissible value for the modulus of convexity $\delta_{d_p}$ is
$$
\delta_{d_p}(\epsilon)=1-\left[1 - \left(\frac{\epsilon}{2}\right)^r\right]^{1/r}=\frac{1}{r2^r}\epsilon^r+o(\epsilon^{2r})\ge k\epsilon^r,
$$
which is increasing in $(0,2]$.
\end{coro}
\begin{proof}
From the previous inequalities follows that if $d_p(\gamma_0,\gamma_1) > \epsilon \max\{d_p(\gamma_0,a),d_p(\gamma_1, a)\}$, then
$$
d_p(a,\gamma_{1/2})^r< \left[1-\left(\frac{\epsilon}{2}\right)^r \right] \max\{d_p(\gamma_0,a),d_p(\gamma_1, a)\}^r.\qedhere
$$
\end{proof}

Note that the formula for the modulus matches those obtained by Clarkson for the $L^p$ measure spaces in \cite{clarkson}.

\subsubsection{Further Inequalities on $\sp$}

\begin{teo}
Let $a\in \sp$, $\gamma:[0,1]\rightarrow \sp$ a geodesic,
$1<p<\infty$. Then there exists a
positive constant $b_p$ such that if $t\in [0,1]$, then
$$
d_p(a,\gamma_t)^r\leq
(1-t)d_p(a,\gamma_0)^r+td_p(a,\gamma_1)^r-w_r(t)b_pd_p(\gamma_0,\gamma_1)^r.
$$
Here $w_r(t)=t^r(1-t)+t(1-t)^r$.
\end{teo}
\begin{proof}
It suffices to prove the assertion for $t\in (0,1)$. Let $f(t)=[\frac{d_p(a,\gamma_{t})}{d_p(\gamma_0,\gamma_1)}]^r$ on
$[0,1]$ and
$$
h(t)=\frac{(1-t)f(0)+tf(1)-f(t)}{w_r(t)}.
$$
Then $h$ is nonnegative in $(0,1)$, since $w_r(t)>0$ in $(0,1)$, and $f$ is convex due to Corollary \ref{puntocurva}. It suffices to prove that the non negative number $b_p:=\inf\{h(t):t\in (0,1)\}$ is strictly positive. Consider
$$\overline{h}(t)=
\left\{%
\begin{array}{lll}
h(0^+) & \mbox{if }&t=0 \\
h(t) & \mbox{if }&0<t<1
\\
h(1^-) & \mbox{if }&t=1 \\
\end{array}%
\right.  , $$
where
$$
h(0^+)=\lim_{t\to 0^+} h(t)=-f(0)+f(1)-f'(0^+),
$$
and
$$
h(1^-)=\lim_{t\to 1^-} h(t)=f(0)-f(1)+f'(1^-).
$$
Then $\overline{h}$ is continuous in $[0,1]$ and it attains its
minimum value $b_p=\overline{h}(c)$ for some $c\in [0,1]$. We claim that $b_p$ is strictly positive.

First assume that $c\in (0,1)$. If $h(c)=0$, then $f(c)=cf(1)+(1-c)f(0)$, and since $f$ is convex and  differentiable, it must be $f(t)=(1-t)f(0)+tf(1)$ for any $t\in (0,1)$, which it is not possible by Theorem \ref{conveunmedio}. Now assume that $c=0$ or $c=1$, and $b_p=h(0^+)=0$ or $b_p=h(1^-)=0$. With an analogous argument one obtains that $f$ is a linear map, which again contradicts Theorem \ref{conveunmedio}.
\end{proof}

\subsection{Projection to convex closed sets}

Let us discuss the properties of the best approximation in convex sets in $\sp$. We start with a brief discussion at the tangent level.

\subsubsection{Linear $p$-orthogonality.}

Let $K\subset \p$ be a convex exponential set. Let $H$ be the tangent space of $K$ at $a=1$, namely $K=\exp(H)$ with $H$ a LTS (see Section \ref{convexos}). Then $a^{1/2}{H} a^{1/2}$ is the tangent space of $K$ at $a\in K$. If $p=2$, then $H$ is clearly complemented. Let us show how to construct a nonlinear supplement for $K$ when $1<p<\infty$.

Since $\lp_h$ is uniformly convex for $1<p<\infty$, there exists \cite{reich}, for any closed convex (in the standard, linear sense) set $C$ in $\lp_h$, a (possibly nonlinear) continuous projection $P_C:\lp_h\to C$ such that
$$
\|x-P_C(x)\|_p\le \|x-y\|_p  \quad \mbox{for any }y\in C,
$$
called the \textit{nearest point projection}. Note that if $C$ is a linear space and $y\in C$, then $y'=P_C(x)+y\in C$ hence
$$
\|x-P_C(x)-y\|_p=\|x-y'\|_p\ge \|x-P_C(x)\|_p,
$$
showing that $P_C\circ (1-P_C)=0$. Clearly $(1-P_C)\circ P_C=P_C-P^2_C=0$ also, hence $P_C$ shares many nice properties with the linear orthogonal projection corresponding to a norm derived from an inner product (which corresponds to $p=2$ for us).

\begin{defi}
Let $1<p<\infty$. Let $K\subset \sp$ be a convex exponential set. Then ${H}_a=a^{1/2}{H}a^{1/2}$ is the tangent space of $K$ at $a\in K$ if $H=T_1K$. Let $P$ be the projection to ${H}$, and let
$$
{H}^{\perp}=H_{1,p}^{\perp}=\{v\in \lp_h: P(v)=0\}
$$
be the Birkhoff orthogonal to $H$. Note that it may contain unbounded self-adjoint operators, i.e. elements of $L^p_h$.
Then
$$
{H}_{a,p}^{\perp}=a^{1/2}H^{\perp}a^{1/2}=\{v\in \lp_h: P(a^{-1/2}{v}a^{-1/2})=0\}
$$
is the Birkhoff orthogonal of ${H}_a$ in $L^p(\M,\tau)_h$. Any element $x\in L^p(\M,\tau)_h$ can be uniquely decomposed as
$$
x=x-P_a(x)+P_a(x)=x_{a,p}^{\perp}+x_{a,p}
$$
where $x_{a,p}=P_a(x)$ is the projection to the completion of ${H}_a$ relative to $d_p$, and $x_{a,p}^{\perp}\in {H}_{a,p}^{\perp}$.
\end{defi}

\begin{lem}
Let $K=\exp(H)\subset \sp$ be a convex exponential set. Then $v\in {H}_{a,p}^{\perp}$ if and only if
$$
\tau\left[ | a^{-1/2} v a^{-1/2} |^{p-1}u^* w  \right]=0\quad \mbox{for any }w\in {H},
$$
where $a^{-1/2} v a^{-1/2}=u |a^{-1/2} v a^{-1/2}|$ is the polar decomposition of $a^{-1/2} v a^{-1/2}$.
\end{lem}
\begin{proof}
It suffices to prove the assertion for $a=1$ and $v\notin H$. For given $w\in H$, consider  $f(t)=\|v-tw\|_p$. The function $f$ is convex and its derivative can be computed with the chain rule and using a supporting functional for $v$ (see, for instance \cite[Theorem 2.3]{abat}),
$$
\dot{f}(0)=\tau\left( \frac{|v|^{p-1}u^*w}{\|v\|_p^{p-1}} \right)
$$
where $v=u|v|$ is the polar decomposition of $v$. Then $t=0$ is a minimum of $f$ if and only if $\tau\left( |v|^{p-1}u^*w\right)=0$.

Now if $\tau\left( |v|^{p-1}u^*w\right)=0$ for any $w\in H$, then $\|v-w\|_p=f(1)\ge f(0)=\|v\|_p$ for any $w\in H$, hence $v\in H^{\perp}$. Likewise, if $v\in H^{\perp}$, for given $w\in H$,  $\|v-tw\|_p\ge \|v\|_p$ for any $t\in\mathbb R$, which says that $f$ has a minimum at $t=0$, so it must be $\tau\left( |v|^{p-1}u^*w\right)=0$.
\end{proof}

\subsubsection{Nonlinear minimization}

Let us consider now the nonlinear projection in the manifold of positive operators (here  $1<p<\infty$). Given the distance $d_p$ in $\sp$, the distance from $a\in \sp$ to a subset $K\subseteq \sp$ is defined according to $d_p(a,K)={\rm{inf}}\{d_p(a,b):b\in K\}$.

\begin{teo}
Let $K= \exp(H)\subset \sp$ be a convex exponential set. Let $b\in \sp$, $b\notin K$. Then $a\in K$ is the best approximation to $b$ in $K$ relative to $d_p$ if and only if the short geodesic $\alpha$ joining $a$ to $b$ in $\sp$ has initial speed $\dot{\alpha}(0)\in {H}_{a,p}^{\perp}$.
\end{teo}
\begin{proof}
Let $c\in K$, let $\beta$ be the short geodesic joining $a$ to $c$ in $K$ ($K$ is convex). Then $f(t)=d_p(b,\beta(t))$ has a global minimum at $t=0$ if and only if $a$ is a best approximation to $b$ in $K$. Now $f$ is convex by Corollary \ref{puntocurva}, hence $t=0$ is minimum of $f$ if and only if $\dot{f}(0)=0$. We may assume (by the invariance of the metric) that $a=1$. Let $c=e^x$ with $x\in {H}$, let $b=e^y$. Then $\beta(t)=e^{t x}$, and $f(t)=\| \ln (e^{-y/2}e^{tx}e^{-y/2}) \|_p=\|\alpha(t)\|_p$, where $\alpha(t)=\ln (e^{-y/2}e^{tx}e^{-y/2})$. Hence, since $\alpha(0)=-y$,
$\dot{f}(0)=\|y\|_p^{1-p}\tau\left( |y|^{p-1}u^* \dot{\alpha}(0)\right)$. Now $e^{\alpha}=e^{-y/2}e^{tx}e^{-y/2}$, hence $(\exp)_{*\alpha}(\dot{\alpha})=e^{\alpha}\,e^{y/2}xe^{-y/2}$, so
$$
e^{y/2}x e^{-y/2}=\int_0^1 e^{-s\alpha}\dot{\alpha}e^{s\alpha}ds
$$
by Lemma \ref{expo}. In particular  $(t=0)$, $e^{y/2}x e^{-y/2}=\int_0^1 e^{sy}\dot{\alpha}(0)e^{-sy}ds$. Since $y$ is self-adjoint, it commutes with $|y|$ and $u$, and also $u=u^*$. Then
$$
\tau(|y|^{p-1}u x)=\int_0^1 \tau(e^{s y}|y|^{p-1}u\dot{\alpha}(0)e^{-s y}) ds=\|y\|_p^{p-1}\dot{f}(0).\qedhere
$$
\end{proof}

The uniqueness of the best approximation $P_K(a)\in K$, and the continuity of the map $P_K$, can be easily derived adapting the proofs of \cite[Theorems 3.15 and 3.17]{conde1} to our context.

\begin{teo}[Best Approximation]

Let $K \subset \sp$ be a convex set, $1< p <\infty$ and $a\in \sp$. Then the best approximation problem has a unique solution in the completion of $K$. In other words,  there is a unique $a_0$ in the completion of $(K,d_p)$ such that  $d_p(a,a_0)=d_p(a, K)$.
\end{teo}
\begin{proof}
Let $\{q_n\}_{n \in \mathbb{N}}$ be a sequence in $K$, such that $d_p(a, q_n)\rightarrow d_p(a, K)$, by Theorem \ref{conveunmedio} we immediately derive that
\begin{equation}\label{centros}
 \frac{1}{2^r}d_p(q_n,q_m)^r\leq  \frac{1}{2}(d_p(q_n,a)^r+d_p(a,q_m)^r)- d_p(a,K)^r,
\end{equation}
where $q_{n,m}=\gamma_{1/2}\in K$ with $\gamma_t$ the geodesic joining $q_n$ and $q_m$.

This implies that $\{q_n\}_{n \in \mathbb{N}}$ is a Cauchy sequence in $K$, hence convergent to some $a_0$ in its completion. By the continuity of the distance we have
$$
d_p(a_0,a)=\lim\limits_{n\to \infty}d_p(q_n, a)=d_p(a,K).
$$
For the uniqueness let $b,c$ in the completion of $K$ with $d_p(b, a)=d=d_p(c,a)$, with $d=d_p(a,K)$. Let $\{b_n\},\{c_n\}$ be Cauchy sequences in $K$ converging to $b,c$ respectively. Replacing $q_n$ and $q_m$ by $b_n$ and $c_n$ respectively in (\ref{centros}) we obtain
$$
d^r\le d_p(a,a_n)^r  \le \frac{1}{2}(d_p(b_n,a)^r+d_p(a,c_n)^r)-\frac{1}{2^r}d_p(b_n,c_n)^r,
$$
where $a_n\in K$ is the midpoint of the geodesic of $K$ joining $b_n$ to $c_n$. Hence
$$
\frac{1}{2^r}d_p(b_n,c_n)^r\le \frac{1}{2}(d_p(b_n,a)^r+d_p(a,c_n)^r)-d^r.
$$
Letting $n\to \infty$ shows that $d_p(b_n,c_n)\to 0$, hence $b=c$.
\end{proof}

\begin{teo}
Let $K \subset \sp$ be a convex set, which is complete for the geodesic distance. Let $1< p <\infty$ and $P_K:\sp\rightarrow K$ the projection onto $K$. Then $P_K$ is continuous.
\end{teo}
\begin{proof}
Let the sequence ${x_n}$ converge to $x$ in $\sp$. Denote $u_n=P_K(x_n)$, which we claim is a Cauchy sequence in $K$. If not, there exist a positive number $\epsilon$ and subsequences ${u_{n_{k}}}$ and ${u_{m_{k}}}$ such that $n_{k}<m_{k}$ and $d_p({u_{n_{k}}},{u_{m_{k}}})\geq \epsilon$ for all $k$. Put $a_k={u_{n_{k}}}, b_k={u_{m_{k}}}$ and $M_k=\max\{d_p(x,a_k),d_p(x,b_k)\}$.
Note that $M_k \rightarrow d_p(x, K)$ as $k\rightarrow \infty$. Now $d_p(x,a_k)\leq M_k$, $d_p(x,b_k)\leq M_k$ and  $d_p({a_{k}},{b_{k}})\geq \frac{\epsilon}{M_k}M_k$. Then, if $m_k\in K$ denotes the midpoint between $a_k,b_k\in\sp$, by Corollary \ref{bco}
$$
d_p(x,K)\le d_p(x, m_k)\leq M_k (1-\delta_{d_p} (\frac{\epsilon}{M_k})).
$$
Hence $\delta_{d_p}(\frac{\epsilon}{M_k})\leq 1-\frac{d_p(x,K)}{M_k}$. Letting $k\rightarrow \infty$, we obtain $\delta(\epsilon)\le 0$, and $\epsilon $ can not be positive.
Thus $\{P_K(x_n)\}$ is a Cauchy sequence in $K$ and therefore converges to a point $z$ in $K$. Since $d_p(x,z)=d_p(x,K)$, then $z=P_K(x)$.
\end{proof}

\begin{rem}
With little effort, many of the results on this paper (for instance, mimimality of the curves $\gamma_{a,b}$ or the convexity of the geodesic distance) can be extended to any tracial gauge norm on the finite von Neumann algebra $\cal M$.
\end{rem}

\bigskip
\bigskip

\noindent
Cristian Conde and Gabriel Larotonda\\
Instituto de Ciencias \\
Universidad Nacional de General Sarmiento \\
J. M. Gutierrez 1150 \\
(B1613GSX) Los Polvorines \\
Buenos Aires, Argentina  \\
and\\
Instituto Argentino de Matem\'atica-CONICET\\
e-mails: cconde@ungs.edu.ar, glaroton@ungs.edu.ar

\end{document}